\documentclass[a4paper,11pt]{article}  
\usepackage{array}                        % array extensions
\usepackage{theorem}                      % theorem extensions
\usepackage{amsmath,amscd,amssymb}        % AMS-LaTeX with CD and symbols
\usepackage[german,english]{babel}        % German & English (default)
 \usepackage[all]{xy}

\newtheorem{lemma}{Lemma}
\newtheorem{proposition}{Proposition}

\newcommand{\CC}{{\mathbb{C}}}

\newcommand{\HH}{{\mathbb{H}}}

\newcommand{\PP}{{\mathbb{P}}}
\newcommand{\QQ}{{\mathbb{Q}}}
\newcommand{\ZZ}{{\mathbb{Z}}}

\newcommand{\RR}{{\mathbb{R}}}
\newcommand{\Eins}{{\mathbf{1}}}

\newcommand{\transp}[1]{{\vphantom{#1}}^t{#1}}

\newcommand{\on}[1]{\operatorname{#1}}
\newcommand{\Gp}{{\widetilde{\Gamma}^{\circ}_{1,p}}}
\newcommand{\Gpe}{{\widetilde{\Gamma}^{\ast}_{1,p}}}
\newcommand{\oGp}{{\Gamma^{\circ}_{1,p}}}
\newcommand{\oGpe}{{\Gamma^{\ast}_{1,p}}}

\newcommand{\Gplev}{{\widetilde{\Gamma}_{1,p}^{\circ}{\scriptstyle (2)}}}
\newcommand{\Gpleverw}{{\widetilde{\Gamma}_{1,p}^{\ast}{\scriptstyle (2)}}}
\newcommand{\oGpleverw}{{\Gamma_{1,p}^{\ast}{\scriptstyle (2)}}}
\newcommand{\oGplev}{{\Gamma_{1,p}^{\circ}{\scriptstyle (2)}}}

\newcommand{\Ap}{ {\cal A}_{1,p}^{\circ}}

\newcommand{\Aplev}{{\cal A}_{1,p}^{\circ}{\scriptstyle (2)}}

\newcommand{\spmat}[1]{\left(\begin{smallmatrix}#1\end{smallmatrix}\right)}

\newcommand{\Mod}[1]{(\on{mod} #1)}
\newcommand{\Sp}{\on{Sp}}

\newcommand{\GL}{\on{GL}}
\newcommand{\SL}{\on{SL}}

\newcommand{\SP}{\on{Sp}}

\newcommand{\Mat}{\on{Mat}}

\newcommand{\tensor} {\otimes}

\newcommand{\smallind}[2]{{#1}_{\scriptscriptstyle{#2}}}
\newcommand{\smallindup}[3]{{#1}_{\scriptscriptstyle{#2}}^{\scriptscriptstyle{#3}}}

\newenvironment{Proof*}[1]{\begin{ProofwCaption}{{#1}}}{\end{ProofwCaption}}
\newenvironment{ProofwCaption}[1]%
  {\addvspace\theorempreskipamount \noindent{\it #1.}\rm}%
  {\newline \mbox{} \hfill \qed \par \addvspace\theorempostskipamount}
\newcommand{\qedsymbol}{\mbox{$\Box$}}
\newcommand{\qed}{ \quad \qedsymbol}
\begin{document}
\selectlanguage{english}
\title{ A discrete extension of $\oGplev$ in $\SP(4,\RR)$ and the
modular form of the Barth-Nieto quintic}
\author{M.~Friedland}
\date{Hannover,13th November 2002}
\maketitle
\begin{quotation}
\begin{center}
{\bf Abstract}
\end{center}
\noindent In this note we construct a maximal discrete extension of
$\Gamma^{\circ}_{1,p}(2)$, the paramodular group with a full level-$2$
structure. The corresponding Siegel variety parametrizes (birationally)
the space of Kummer surfaces associated to $(1,p)$-polarized abelian
surfaces with a level-$2$ structure. In the case $p=3$ this is related
to the Barth-Nieto quintic and in this case we also determine the space
of cusp forms of weight 3.
\end{quotation}
\section{Introduction}
\label{intro}
Let $p \ge 3$ prime. The {\it paramodular group}
$\oGp$ is defined as the subgroup
\[
\oGp = \left\{ g \in \SP(4,\QQ);~g \in \spmat{
\ZZ & \ZZ & \ZZ & p\ZZ \cr
p\ZZ & \ZZ & p\ZZ & p\ZZ \cr
\ZZ & \ZZ & \ZZ & p\ZZ \cr
\ZZ & p^{-1}\ZZ & \ZZ & \ZZ} 
\right\}
\] of $\SP(4,\QQ)$.
This group acts on the Siegel upper halfspace
\[
\HH_2 := \left\{ \tau \in \Mat(2,\CC);~\tau = \transp{\tau},~\on{Im}
    \tau > 0 \right\}
\]
by
\[
\spmat{A & B \cr C & D} :\left\{  \begin{array}{rcl} \HH_2
    &\rightarrow & \HH_2 \cr
    \tau &\mapsto & (A\tau + B)(C\tau + D)^{-1}
    \end{array} \right.
\] (where $A,B,C,D$ are $2\times 2$-blocks).
With this action, $\Ap :=
\oGp\backslash \HH_2$ is the moduli space of $(1,p)$-polarized abelian surfaces.  
Likewise, one obtains the moduli space of $(1,p)$-polarized abelian surfaces
with level 2 structure $\Aplev$ by dividing $\HH_2$ by the action of 
\[
\oGplev = \left\{ g\in \oGp; ~g - \Eins_4 \in \spmat{2\ZZ & 2\ZZ & 2\ZZ & 2p\ZZ \cr
2p\ZZ & 2\ZZ & 2p\ZZ & 2p\ZZ \cr
2\ZZ & 2\ZZ & 2\ZZ & 2p\ZZ \cr
2\ZZ & 2p^{-1}\ZZ & 2\ZZ & 2\ZZ} \right\}.
\]
The group $\oGp$ is conjugate via $R_p := \on{diag}(1,1,1,p)$ to the
symplectic group 
\[
\Gp = \SP(\Lambda_p,\ZZ) := \left\{g \in \GL(4,\ZZ)~;~g\Lambda_p \transp{g}
  = \Lambda_p \right\},
\] where $\Lambda_p$ is the symplectic form 
\[
\Lambda_p = \spmat{
  0&0&1&0\cr
  0&0&0&p\cr
  -1&0&0&0\cr
  0&-p&0&0}.
\]
Under this isomorphism, $\oGplev$ is identified with the group
$\Gplev$ consisting of all elements $g$ in $\SP(\Lambda_p,\ZZ)$ with $g\equiv \Eins_4~\Mod{2}$. 
The case $p=3$ is of special interest. Barth and Nieto showed in \cite{BN} that the quintic 
\[
        N = \left\{ \sum\limits^5_{i=0} u_i = \sum\limits^5_{i=0} \frac{1}{u_i} 
= 0\right\} \subset \PP^5
\] 
parametrizes birationally the space of 
Kummer surfaces associated to abelian surfaces with
$(1,3)$-polarization and a level 2 structure. Moreover, $N$
has a smooth model which is Calabi-Yau. From this, Barth and Nieto
deduced that the space ${\cal A}_{1,3}^{\circ}(2)$ also has
a smooth model that is Calabi-Yau. So one may ask to determine the (up
to a scalar) weight 3 cusp form with respect to the modular group
$\oGplev$. This was done in \cite{GH1}. The cusp form in question was
shown to be $\smallindup{\Delta}{1}{3}$, where $\smallind{\Delta}{1}$ is a
cusp form of weight 1 with respect to the paramodular group
$\Gamma^{\circ}_{1,3}$ with a character of order
6.\\
In \cite{GH2} Gritsenko and Hulek showed that the Kummer surfaces which are associated to a $(1,p)$-polarized
abelian and to its dual are isomophic. This turns our attention to the
{\it Fricke-involution}, which extends $\oGp$ to $\oGpe$ and
identifies a polarized abelian surface with its dual. So it is a natural
question to ask if $\oGplev$ can be extended uniquely to a group $\oGpleverw$ in such a way 
that the diagram 
$$
\xymatrix{
\oGplev~  \ar@{^{(}->}[r]\ar@{^{(}->}[d] & ~\oGpleverw \ar@{^{(}->}[d] \cr
    \oGp~ \ar@{^{(}->}[r] & ~\oGpe
}
$$
commutes. If so, is $\smallindup{\Delta}{1}{3}$ still a cusp form 
with respect to ${\Gamma_{1,3}^{\ast}({\scriptstyle 2})}$?
We will give answers to this question in this note.
\section{The maximal discrete extension $\oGpe$ of $\oGp$ in
  $\Sp(4,\RR)$}  
\label{sec:2}
A  maximal discrete extension $\oGpe$ of $\oGp$ in
$\SP(4,\RR)$ is defined in \cite{GH2}. This group still acts on
$\HH_2$ and the quotient $\oGpe \backslash \HH_2$ has a
moduli theoretic meaning: $\oGpe \backslash \HH_2$ is birationally the moduli
space of Kummer surfaces associated to abelian surfaces with a
$(1,p)$-polarisation (for details see \cite{GH2}).\\
We will construct an extension of $\Gplev$ in a natural way. For this,
we summarize the construction of $\oGpe$. Let $x,y \in \ZZ$ with $xp - y = 1$ and consider
the matrix
\[
        \widehat{V}_{p} = \spmat{        px & -1 & 0 & 0 \cr
                -yp & p & 0 & 0 \cr
                0   & 0 & p & yp \cr
                0 & 0 & 1 & px}.
\]
Let
\[
        V_p = \frac{1}{\sqrt{p}}\widehat{V}_{p} \in \SP(4,\RR).
\]
Then it is easy to see that  $V_p^2 \in \oGp$ and $V_p \oGp V_p = \oGp$.
So the matrix $V_p$ defines an involution modulo $\oGp$  and $\oGpe :=
\left< \oGp,V_p\right>$\label{defoGp} is a normal extension of $\oGp$
with index 2. By \cite{K} this is the only non-trivial discrete extension of $\oGp$ in
$\SP(4,\RR)$.\\
With
\[\bar{V}_p = \frac{1}{\sqrt{p}}\spmat{0 & 1 & 0 & 0 \cr   
                p & 0 & 0 & 0 \cr
                0 & 0 & 0 & p \cr
                0 & 0 & 1 & 0 }, \label{defVp}
\]
the coset $V_p \oGp$ can also be written as 
\[
        V_p \oGp =    \bar{V}_p  \oGp.
\]
To understand how $\bar{V}_p$ acts on $\Aplev$ let $E =
\on{diag}(1,p)$ and $\tau =
\spmat{\tau_1 & \tau_2 \cr \tau_2 & \tau_3} \in \HH_2$ be a point
corresponding to the $(1,p)$-polarized abelian surface $X = \CC^2/L$,
where the lattice  $L$ is given by the normalized period matrix
$\Omega = (E,\tau)$ and the hermitian form $H$, defining the
polarization of $X$, is given by $(\on{Im}
\tau)^{-1}$ with respect to the standard basis of $\CC^2$. The
polarization  $H$ defines an isogeny
\[
        \lambda_{H} : \left\{ \begin{array}{rcl} X & \rightarrow & \widehat{A} =
\on{Pic}^{\circ}A \cr
                        x & \mapsto & T^{\ast}_x {\cal L} \tensor {\cal L}^{-1}
                        \end{array}
                        \right.
\]
where ${\cal L}$ is a line bundle, which represents the polarization
$H$ and $T_x$ is the translation by $x$. The map $\lambda_H$ depends
only on the polarization, not on the choice of the line bundle ${\cal
  L}$. The kernel $\on{ker}\lambda_H$ is (non-canonically)
isomorphic to $\ZZ_p \times \ZZ_p$, so this defines a quotient map \[
        \lambda_p : X \rightarrow X/\on{ker} \lambda_H = \widehat{X},
\]
where $\widehat{X}$ is the dual abelian surface of $X$, which corresponds
to the period matrix
\[
        \Omega' = \spmat{ p & 0 & p\tau_1 & \tau_2\cr
                         0 & 1 & \tau_2 & \tau_3/p }.
\]
The identity     
\[
        \spmat{ 0 & 1 \cr 1 & 0} \spmat{ p & 0 & p\tau_1 & \tau_2\cr
                         0 & 1 & \tau_2 & \tau_3/p } \left({ {\spmat{0 & 1
\cr 1 & 0}} \atop {0_2}} {{0_2} \atop {\spmat{0 & 1 \cr 1 & 0}}} \right) =
\left(\spmat{ 1 & 0 \cr 0 & p}, \bar{V}_p(\tau) \right)
\]
shows that $\widehat{X}$ is $(1,p)$-polarized and hence the action of $\bar{V}_p$ induces a morphism
\[
        \varphi(p) : \left\{\begin{array}{rcl} 
                        \Ap & \rightarrow & \Ap \cr
                        (X,H) & \mapsto & (\widehat{X},\widehat{H})
                                \end{array}
                        \right.
\]
which maps an abelian surface to its dual.\\
Since it is easier to work with matrices with entries in $\ZZ$, let us
consider
\[
        \widetilde{W}_p =  R_p \bar{V}_p R_p^{-1} = \frac{1}{\sqrt{p}}\spmat{0 & 1 & 0 & 0 \cr   
                p & 0 & 0 & 0 \cr
                0 & 0 & 0 & 1 \cr
                0 & 0 & p & 0 }.\label{defWp} 
 \]
We define $\Gpe := \left< \Gp, \widetilde{W}_p \right>\label{defGpe}$. 
Recall that $\Gplev$ is the kernel of the surjection
$\pi:~ \Gp \rightarrow \Sp(4,\ZZ_2) \simeq S_6,~ \pi(M) =
\overline{M}$. (Here and henceforth we write $\overline{M}$ for reduction modulo 2 of
an integer-valued matrix $M$.) 
\begin{lemma}
Let 
\[      \iota = \spmat{\spmat{0 & 1\cr 1 & 0} & 0_2 \cr
                        0_2 & \spmat{0 & 1\cr 1 & 0}}. 
\] 
The map  
\[
        \pi^{\ast} : \left\{ \begin{array}{rcl}
                                \Gpe & \rightarrow & \SP(4,\ZZ_2) \cr
                                g & \mapsto & \left\{ \begin{array}{c} \pi(g) \vspace{6pt}\cr
                                                \pi(g\cdot\widetilde{W}_p)\cdot\iota
                                                \end{array} 
                                                                ~\mbox{if}~ 
                                                \begin{array}{l} g \in \Gp \vspace{6pt}\cr 
                                                                g \in \Gpe\backslash \Gp 
                                                \end{array}
                                                \right.
                        \end{array}
                        \right.\label{defpiast}
\]
is a homomorphism which extends the map $\pi$.
\end{lemma}
\begin{Proof*}{Proof}
It is easy to see that the equation 
\begin{center}
         $\iota \cdot \pi^{\ast}(\widetilde{W}_p \cdot g \cdot \widetilde{W}_p) \cdot 
\iota = \pi^{\ast}( g )~~\mbox{ for all } g \in \Gp ~~~(\ast)$
\end{center}
holds. Namely, let $g = \spmat{A & B \cr C & D} \in \Gp,~ A = (a_{ij})_{1\leq 
i,j \leq 2},\dots, D = (d_{ij})_{1\leq 
i,j \leq 2}$. By \cite[Proposition I.1.16]{HKW} we have $a_{21}\equiv
b_{21} \equiv c_{21} \equiv d_{21} \equiv 0~ \Mod{p})$ and 
\[
        \begin{array}{rcl}
         \iota\cdot \pi^{\ast}(\widetilde{W}_p \cdot g \cdot
 \widetilde{W}_p)\cdot \iota &=&
 \iota \cdot \pi(\widetilde{W}_p \cdot g \cdot \widetilde{W}_p)\cdot \iota\cr
                 &=& \iota \cdot \pi \left(\spmat {\frac{1}{p}\spmat{0 & 1 \cr p & 0}\cdot A \cdot
 \spmat{0 & 1 \cr p & 0} & \frac{1}{p}\spmat{0 & 1 \cr p & 0}\cdot B \cdot
 \spmat{0 & 1 \cr p & 0}\vspace{5pt}\cr \frac{1}{p}\spmat{0 & 1 \cr p & 0}\cdot C \cdot
 \spmat{0 & 1 \cr p & 0} &
         \frac{1}{p}\spmat{0 & 1 \cr p & 0}\cdot D \cdot \spmat{0 & 1 \cr p &
 0}}\right)\cdot \iota  \vspace{5pt} \cr
                &=&\iota \cdot \spmat {\overline{\spmat{a_{22} & a_{21}/p \cr pa_{12} &
a_{11}}} & \overline{\spmat{b_{22} & b_{21}/p \cr pb_{12} & b_{11}}}\vspace{5pt}\cr \overline{\spmat{c_{22} & c_{21}/p \cr pc_{12} & c_{11}}} &
        \overline{\spmat{d_{22} & d_{21}/p \cr pd_{12} & d_{11}}}
}       \cdot \iota \vspace{5pt} \cr
                &=&  \spmat {\overline{
\spmat{ \vphantom{b}a_{11} & a_{12} \cr \vphantom{b} a_{21} &
a_{22}}}  & \overline{\spmat{b_{11} & b_{12} \cr
b_{21} & b_{22}}} \vspace{5pt}\cr \overline{\spmat{\vphantom{b}c_{11} & c_{12}\cr \vphantom{b}c_{21} & c_{22}}}  &
         \overline{\spmat{d_{11} & d_{12} \cr
d_{21} & d_{22}}}  }\vspace{5pt} \cr
                &=& \pi^{\ast}(g).
        \end{array}
\]
Since $\Gp$ is normal in $\Gpe$ and $\widetilde{W}_p^2 = \Eins_4$,
$(\ast)$ is equivalent to 
\[
        \iota\cdot \pi^{\ast} (\widetilde{W}_p \cdot h) = \pi^{\ast} (h\cdot
\widetilde{W}_p)\cdot \iota ~~\mbox{ for all } h \in \Gpe\backslash \Gp. ~~~(\ast \ast)
\]
Let $g \in \Gp,~h_1,h_2 \in \Gpe \backslash \Gp$. Then 
\[ \begin{array}{rrcl}
                (i) & \pi^{\ast}(g\cdot h_1) & = & \pi(g\cdot h_1 \cdot
\widetilde{W}_p) \cdot \iota \cr
                && = &\pi(g) \cdot \pi(h_1\cdot \widetilde{W}_p) \cdot \iota \cr
                && = &\pi^{\ast}(g) \cdot \pi^{\ast}(h_1)\cr
                (ii) & \pi^{\ast}(h_1\cdot g) & = & \pi(h_1\cdot g \cdot
\widetilde{W}_p) \cdot \iota \cr
                && = &\pi(h_1\cdot \widetilde{W}_p \cdot \widetilde{W}_p \cdot
g\cdot \widetilde{W}_p) \cdot \iota \cr
                && = & \pi(h_1\cdot\widetilde{W}_p) \cdot \pi (\widetilde{W}_p
\cdot g \cdot \widetilde{W}_p) \cdot \iota \cr
                &&= & \pi^{\ast}(h_1) \cdot \iota \cdot \pi(\widetilde{W}_p \cdot g
\cdot \widetilde{W}_p) \cdot \iota \cr
                &&  \stackrel{(\ast)}{=} & \pi^{\ast}(h_1) \cdot
                \pi^{\ast}(g)\cr
(iii)& \pi^{\ast}(h_1 \cdot h_2) &=& \pi^{\ast}(h_1 \cdot
\widetilde{W}_p \cdot\widetilde{W}_p \cdot h_2)\cr
                && = &\pi^{\ast}(h_1 \cdot \widetilde{W}_p) \cdot
\pi^{\ast}(\widetilde{W}_p \cdot h_2) \cr
                && =& \pi^{\ast}(h_1) \cdot \iota \cdot
\pi^{\ast}(\widetilde{W}_p \cdot h_2) \cr
                &&  \stackrel{(\ast\ast)}{=}& \pi^{\ast}(h_1) \cdot \pi^{\ast}(h_2 \cdot
\widetilde{W}_p) \cdot \iota \cr
                && =& \pi^{\ast}(h_1)\cdot \pi^{\ast}(h_2).     
        \end{array}
\]
\end{Proof*}
\begin{proposition}
There is exactly one group $\Gpleverw$ such that the diagram 
$$
\xymatrix{
& 1 \ar@{->}[d]&1 \ar@{->}[d] &  &\cr
1 \ar@{->}[r] & \Gplev \ar@{->}[d] \ar@{->}[r] & \Gpleverw \ar@{->}[d]\ar@{->}[r] & 
\ZZ_2 \ar@{=}[d]\ar@{->}[r] & 1\cr
1 \ar@{->}[r] & \Gp \ar@{->}[d]\ar@{->}[r] & \Gpe \ar@{->}[d]\ar@{->}[r] & 
\ZZ_2  \ar@{->}[r] & 1\cr
 & \SP(4,\ZZ_2) \ar@{->}[d] \ar@{=}[r] & \SP(4,\ZZ_2) \ar@{->}[d]&  & \cr
& 1&1  &   &
}
$$
commutes with exact rows and columns.
\end{proposition}
\begin{Proof*}{Proof}
Let $\varphi$ be a homomorphism such that the diagram
\[
\xymatrix{
\Gp  \ar@{^{(}->}[r]\ar@{->>}[d]^{\pi} & \Gpe \ar@{->>}[d]^{\varphi} \cr
    \Sp(4,\ZZ_2) \ar@{=}[r] & \SP(4,\ZZ_2)
}
\]
commutes. It is enough to show that one necessarily has $\varphi(\widetilde{W}_p)~=~\pi^{\ast}(\widetilde{W}_p)~=~\iota.$\\
First, $\varphi(\widetilde{W}_p)$ is an involution in $\Sp(4,\ZZ_2)$
since $\widetilde{W}_p^2 = \Eins_4$. Moreover, we have
\[
\begin{array}{lr}
\varphi(g\cdot\widetilde{W}_p\cdot g^{-1}) =
\pi(g)\cdot\varphi(\widetilde{W}_p)\cdot\pi(g)^{-1} & \forall~g \in
\Gp,
\end{array}
\]
so 
\[
        g \in \on{centr}(\widetilde{W}_p,\Gpe) \Rightarrow \pi(g) \in
\on{centr}(\varphi(\widetilde{W}_p),\Sp(4,\ZZ_2))
\]
($\on{centr}(x,G)$ means the centralisor of $x$ in $G$). The matrices
\[
        h_1 = \spmat{    1 & 0 & 0 & 1 \cr
                        0 & 1 & p & 0 \cr
                        0 & 0 & 1 & 0 \cr
                        0 & 0 & 0 & 1 } ~\mbox{ and }~
        h_2= \spmat{     1 & 0 & 0 & 0 \cr
                        0 & 1 & 0 & 0 \cr
                        0 & 1 & 1 & 0 \cr
                        p & 0 & 0 & 1 } 
\] are in  $\on{centr}(\widetilde{W}_p,\Gpe)$. Since  
\[
        \pi(h_1) = \spmat{ 1 & 0 & 0 & 1 \cr
                        0 & 1 & 1 & 0 \cr
                        0 & 0 & 1 & 0 \cr
                        0 & 0 & 0 & 1 }
\]
must lie in $\on{centr}(\varphi(\widetilde{W}_p),\SP(4,\ZZ_2))$, it follows 
(with $\varphi(\widetilde{W}_p) = \spmat{A & B \cr C & D}$), that $C =
0_2$. Likewise, one argues with $\pi(h_2)$ that necessarily $B 
= 0_2$ holds.\\
Since $\varphi(\widetilde{W}_p)$ is an involution in $\Sp(4,\ZZ_2)$,
we deduce that $A$ (and hence $D$) has to be an involution in
$\SL(2,\ZZ_2)$ (i.~e.~equal to $\Eins_2,
\spmat{0&1\cr1&0},\spmat{1&1\cr0&1}$ or $\spmat{1&0\cr1&1} $).\\
Now it is easy to see that necessarily $A=D=\spmat{0&1\cr 1&0}$
holds. To see this first assume $A=\Eins_2$ (and consequently
$D=\Eins_2$). Then necessarily $\pi(\widetilde{W}_p\cdot h
\cdot \widetilde{W}_p) = \pi(h)$ must hold for all $h \in \Gp$. But
this is not the case for
\[
        \spmat{      1 & 1 & 0 & 0 \cr
                        0 & 1 & 0 & 0 \cr
                        0 & 0 & 1 & 0 \cr
                        0 & 0 & -p & 1 } ~\in ~\Gp.
\]    
Similarly one argues with    
\[
        \spmat{  1 & 0 & 0 & 0 \cr
                        p & 1 & 0 & 0 \cr
                        0 & 0 & 1 & -1 \cr
                        0 & 0 & 0 & 1 }  ~\in ~\Gp 
\]
to exclude the cases $A =\spmat{ 1&1\cr0&1}$
and $A=\spmat{1&0\cr1&1}$. This shows that the only possibility is
$\Gpleverw = \ker(\pi^{\ast})\label{defkererw}.$
\end{Proof*}
For future use we note that the kernel $\ker(\pi^{\ast})$ ist generated by $\Gplev$ and
$\widetilde{\kappa}_p$ with $\widetilde{\kappa}_p = \widetilde{W}_p \cdot g$ and 
$g \in \Gp$, $\pi^{\ast}(g) =\iota$. Of course we have $G = \left<
\Gplev,~\widetilde{W}_p \cdot g_1\right>
= \left< \Gplev,~\widetilde{W}_p \cdot g_2\right>$ for $g_1,g_2 \in \Gp$ with $\pi^{\ast}(g_1) =\pi^{\ast}(g_2) =
\iota$, so we can choose
\[
        \begin{array}{rcl}
        \widetilde{\kappa}_p &=& \widetilde{W}_p \cdot \spmat{ p-1 & 2-p & 0 & 0 \cr 
                                                        p & 1-p & 0 & 0 \cr 
                                                        0 & 0 & p-1 & 1 \cr 
                                                        0 & 0 & p(2-p) & 1-p}
\vspace{4pt}\cr
                 &=& \frac{1}{\sqrt{p}} \widetilde{V}_p,~ \widetilde{V}_p = \spmat{      p & 1-p & 0 & 0 \cr 
                                                p(p-1) & p(2-p) & 0 & 0 \cr 
                                                0 & 0 & p(2-p) & 1-p \cr 
                                                0 & 0 & p(p-1) & p}
        \end{array}
\]
as a generator.
\section{ The modular form $\smallindup{\Delta}{1}{3}$}
\label{sec:3}
It was shown in \cite{GN} that 
\[
        \smallind{\Delta}{1}(\tau) = q^{1/6}r^{1/2}s^{1/2} \hspace{-.85cm}\prod_{\stackrel{n \geq 0,m \geq 0, l \in
\ZZ}{ (l < 0 \on{if} n =m=0)}}\hspace{-.85cm}(1-q^nr^ls^{3m})^{f(nm,l)}
\]
with $q = e^{2\pi i \tau_1},~ r = e^{2 \pi \tau_2},~ s = e^{2 \pi i
  \tau_3}$ and
\[
        \begin{array}{l}
        \sum_{n \geq 0,l} f(n,l)q^nr^l = \cr
         r^{-1} \left( \prod_{n \geq 1}(1 +
q^{n-1}r)(1 + q^nr^{-1})(1-q^{2n-1}r^2)(1-q^{2n-1}r^{-2}) \right)^2 
        \end{array}\label{defDelta}
\]
is a cusp form of weight one  with respect to  $\Gamma_{1,3}^{\circ}$
with a character $\smallind{\chi}{6}\label{defChar}$ of order 6. Then
$\smallindup{\Delta}{1}{3}$ is a cusp form with respect to
$\Gamma_{1,3}^{\ast}$ with a character
$\smallind{\chi}{(2,-)}\label{defChar2}$ of order two and
$\smallind{\chi}{(2,-)}(\bar{V}_3) = -1$. The character $\smallind{\chi}{(2,-)|\Gamma_{1,3}^{\circ}}$
arises in the following way: 
$$
\xymatrix{
1 \ar@{->}[r] & \Gamma_{1,3}^{\circ}({\scriptstyle 2}) \ar@{->}[r] &
{\Gamma^{\circ}_{1,3}} \ar@{->}[r] \ar@{->}[rd]^{ \smallind{\chi}{(2,-)}} & \SP(4,\ZZ_2)\simeq S_6 \ar@{->}[r]\ar@{->}[d]^{\on{sgn}} & 1 \cr
 &&&\left\{\pm 1 \right\}&
}
$$
(for all this see  (see \cite{GH2})).
\begin{proposition}
Let $S_3(\Gamma_{1,3}^{\ast}({\scriptstyle 2}))$ be the vectorspace of
cusp forms of weight 3 with respect to $\Gamma_{1,3}^{\ast}({\scriptstyle
2})$. Then 
\[
        S_3(\Gamma_{1,3}^{\ast}({\scriptstyle 2})) = \CC \cdot \Delta^3_1.
\]
\end{proposition}
\begin{Proof*}{Proof}
Consider $\Gamma_{1,3}^{\ast}({\scriptstyle 2})
= \left<\Gamma_{1,3}^{\circ}({\scriptstyle 2}),
R_3^{-1}\cdot\widetilde{\kappa}_3 \cdot R_3\right>$. From the diagram above
we obtain $\smallind{\chi}{(2,-)|\Gamma_{1,3}^{\circ}({\scriptstyle
    2})} \equiv 1$. We have  
$R_3^{-1}\cdot \widetilde{\kappa}_3 \cdot R_3 = \bar{V}_3 \cdot g$
with  \[
        g = \spmat{ 2 & -1 & 0 & 0\cr
                3 & -2 & 0 & 0\cr
                0 & 0 & 2 & 3\cr
                0 & 0 & -1 &-2}
\in \Gamma_{1,3}^{\circ}.
\]
Moreover we have $\pi(R_3 \cdot g \cdot R_3^{-1}) = \iota$ and it is easy to see that
$\on{sgn}(\iota) = -1$, so
\[
        \begin{array}{rcl}
        \smallind{\chi}{(2,-)}(R_3^{-1}\cdot \widetilde{\kappa}_3\cdot R_3) &=&
\smallind{\chi}{(2,-)}(\bar{V}_3)\cdot \smallind{\chi}{(2,-)}(g) \cr
                        &=& -1\cdot \on{sgn}(\iota)\cr
                        &=& 1.
        \end{array}
\]
This shows that $\Delta_1^3$ is a cusp form of weight 3 with respect
to $\Gamma_{1,3}^{\ast}({\scriptstyle 2}))$.
By \cite{GH2} and \cite{BN} the moduli space 
$\Gamma_{1,3}^{\ast}({\scriptstyle 2}) \backslash \HH_2$ is birationally
equivalent to a Calabi-Yau variety. Hence the result follows from
Freitag's extension theorem (cf.~\cite[Hilfssatz 3.2.1]{F}) which says that for any discrete group $\Gamma$ the
space $S_3(\Gamma)$ is isomorphic to $H^{3,0}(\widetilde{{\cal
  A}}(\Gamma))$ for any smooth projective model $\widetilde{{\cal
  A}}(\Gamma)$ of ${\cal A}(\Gamma) =
\Gamma \backslash \HH_2$.
\end{Proof*}

\end{document}